\documentclass[12pt,a4paper]{amsart}
\usepackage[utf8]{inputenc}

\usepackage{graphicx,color}
\usepackage{amsmath}
\usepackage{amssymb}
\usepackage{amsfonts}
\usepackage{amsthm}
\usepackage{amscd}
\usepackage{ae}
\usepackage[T1]{fontenc}

\newcommand{\dist}{\mathrm{dist}}

\renewcommand\arg{\operatorname{arg}}

\renewcommand{\bar}[1]{\overline{#1}}

\newcommand{\diam}{{\rm diam \,}}

\theoremstyle{plain}
\newtheorem{theorem}{Theorem}

\newtheorem{lemma}[theorem]{Lemma}

\theoremstyle{definition}

\author[L. Hitruhin]{Lauri Hitruhin}
\email{lauri.hitruhin@helsinki.fi}

\title[Dimension compression and expansion]{Dimension compression and expansion under homeomorphisms with exponentially integrable distortion}

\begin{document}

\begin{abstract}
We improve both dimension compression and expansion bounds for homeomorphisms with $p$-exponentially integrable distortion. To the first direction  we also  introduce estimates for the compression multifractal spectra, which will be used to estimate compression of dimension, and for the rotational multifractal spectra. For establishing the expansion case we use the multifractal spectra of the inverse mapping and construct examples proving sharpness.
\end{abstract}

\maketitle
\section{Introduction}
Let $f: \mathbb{C} \to \mathbb{C}$ be  a homeomorphism with $p$-exponentially integrable distortion for a given $p>0$. For such mappings the sharp pointwise bounds  for compression and stretching
\begin{equation}\label{PisteVen}
e^{-\frac{\bar{C}}{p}\log^{2}\left( \frac{1}{|z|} \right)} \lesssim |f(z)|  \lesssim \frac{1}{\log^{\frac{p}{2}}\left( \frac{1}{|z|} \right)},
\end{equation}
are established in, respectively, \cite{HK} and \cite{OZ}. Using the lower bound of \eqref{PisteVen} Zapadinskaya in \cite{Zap} established dimension compression estimate for these mappings. She proved that given a sufficiently big set $A$ in the Hausdorff measure sense, that is $$H^s(A)>0$$ for given $s\in (0,2)$, then also the image set is big in the sense that $$H^h(f(A))>0$$
for the gauge function
\begin{equation}\label{GaugeZAP}
h(t)=e^{-cs\sqrt{p} \log^{\frac{1}{2}} \left( \frac{1}{t}  \right)}.
\end{equation}
Furthermore, it was shown by Clop and Herron in \cite{CH}, see also slightly weaker original example by Zapadinskaya in \cite{Zap}, that  given any $s\in (0,2)$
there exists a mapping with $p$-exponentially integrable distortion and a set $A$ such that $H^{s}(A)>0$ but $H^{\bar{h}}(f(A))=0$, where the gauge function is 
\begin{equation}\label{GaugeCH}
h(t)=e^{-s\sqrt{p} \frac{2+\epsilon}{2-s}\log^{\frac{1}{2}}\left( \frac{1}{t} \right)}.
\end{equation}
Thus we see that the compression result by Zapadinskaya is asymptotically of correct order when $s\to 0$, but when $s\to 2$ there is a gap of order $(2-s)^{-1}$ in the exponent. This is a natural consequence of Zapadinskaya using the sharp pointwise contraction bound \eqref{PisteVen}, which intuitively speaking captures the correct behaviour only when the set under study has the Hausorff dimension zero, to study compression in a larger scale. Thus it is evident that in order to obtain better bounds for compression of measure for big sets  one should  study local compression properties instead of using the sharp pointwise bound. Finding the optimal local compression for any given dimension $s\in (0,2)$ is called  establishing the compression multifractal spectra for this class of mappings.

\begin{theorem}\label{VenytysT}
Let $p>0$, $s\in (0,2)$ and $f: \mathbb{C} \to \mathbb{C}$ be a homeomorphism with $p$-exponentially integrable distortion. Assume that for any $z\in A$ there exists a sequence $\lambda_{z,n}$, which satisfy $|\lambda_{z,n}|\to 0$ when $n\to \infty$, such that
\begin{equation}\label{MSVEN}
|f(z+\lambda_{z,n})-f(z)| \leq  e^{-\frac{\bar{C}(2-s)}{2p} \log^2 \left( \frac{1}{|\lambda_{z,n}|} \right)},
\end{equation}
where $\bar{C}$ is the constant from the pointwise bound \eqref{PisteVen}. Then 
\begin{equation}\label{DIAMA}
\dim(A) \leq 2- \frac{\bar{C}\tilde{C} (2-s)}{\pi + 2(e-1)},
\end{equation}
where $\tilde{C}$ is the constant from the modulus inequality \eqref{MODLEMMA}.
\end{theorem}
We believe that the dimension bound \eqref{DIAMA} can be improved to the form $\dim(A) \leq s$, in which case  the optimal pointwise contraction can only occur in sets with zero Hausdorff dimension, as for mappings with $p$-integrable distortion, see \cite{H4}. Furthermore, the bound \eqref{DIAMA} shows that if we want compression to happen in big sets, that is, $\dim(A)$ is close to $2$, then  $s$ must be close to $2$. When this happens the compression in \eqref{MSVEN} gets considerably weaker than in the pointwise case \eqref{PisteVen}. Hence the form given in \eqref{DIAMA} achieves better asymptotical behaviour when $s\to 2$ than the pointwise bound \eqref{PisteVen} and can thus be used to obtain better dimension compression in this situation.

\begin{theorem}\label{DimensioT}
Fix $p>0$, assume that $s$ is close to $2$ and let $f: \mathbb{C} \to \mathbb{C}$ be a homeomorphism with $p$-exponentially integrable distortion.  Define the gauge function
\begin{equation}\label{GaugeAlku}
h(t)= e^{- \frac{C\sqrt{p \tilde{C}}  s}{\sqrt{2-s+ \epsilon}} \log^{\frac{1}{2}}  \left(\frac{1}{t}\right)  },
\end{equation}
where $\epsilon>0$ can be chosen arbitrary small. Then every set $A \subset \mathbb{C}$ for which $$H^{h}(f(A))=0$$
satisfies $$H^{s}(A)=0.$$
\end{theorem}
Comparing the gauge function \eqref{GaugeAlku} to the gauge function \eqref{GaugeCH} for which Clop and Herron constructed examples we see that there is a difference of a square root in $2-s$ term. Thus there is still small gap left between positive result and examples. But on the other hand our result greatly improves, when $s$ is big,  the previous positive result of Zapadinskaya \eqref{GaugeZAP} which does not have term $2-s$ at all. It is an interesting question whether the square root in the term $2-s$ is needed or not? \\
\\
Our method can also be used to study the rotational multifractal spectra, which has earlier been studied for quasiconformal mappings in \cite{AIPS} and for homeomorphisms with $p$-integrable distortion in \cite{H4}.

\begin{theorem}\label{MSPEKTK}
Let $p>0$, $s\in (0,2)$ and  $f: \mathbb{C} \to \mathbb{C}$ be a homeomorphism with $p$-exponentially integrable distortion and fix a branch of the argument. Assume that for every $z\in A$ there exists a sequence $\lambda_{z,n}$, for which $|\lambda_{z,n}| \to 0$ when $n\to \infty$, such that 
\begin{equation}\label{KIERTO}
|\arg(f(z+\lambda_{z,n})-f(z))|\geq  \frac{(2-s)\sqrt{\tilde{C} \pi}\bar{C}}{\sqrt{2} p} \log^2 \left( \frac{1}{|\lambda_{z,n}|} \right).
\end{equation}
 Then we can bound the dimension of the set $A$ by 
\begin{equation}\label{DIAMA2}
\dim(A) \leq 2- \frac{\bar{C}\tilde{C} (2-s)}{\pi +2(e-1)}.
\end{equation}
\end{theorem}
This theorem couples size of a set $A$ and the maximal rotation that can occur locally in it. We believe that constants in the rotation bound \eqref{KIERTO} and the dimension bound \eqref{DIAMA2} can be improved, but noteworthly the sharp constant is not known even in the optimal pointwise bound
 $$|\arg(f(z))| \leq \frac{c}{p} \log^2 \left( \frac{1}{|z|} \right)$$ established in Theorem 1 in \cite{H2}. However, already presented form  \eqref{DIAMA2} gives a similar asymptotical behaviour when $s\to 2$ as the compression multifractal spectra.  \\
\\
As mentioned before Clop and Herron in \cite{CH} and Zapadinskaya in \cite{Zap} have constructed examples for dimension compression, which could also be used to derive examples for the compression multifractal spectra. But since this connection has not been written down, and more importantly as the rotational part has not been covered before, we briefly provide examples towards optimality.

\begin{theorem}\label{Esimerkki}
Fix $p>0$, $s\in (0,2)$ and a branch of the argument. There exists a homeomorphism with $p$-exponentially integrable distortion and a set $A$ with Hausdorff dimension  $s$ such that for each point $z\in A$ there exists a sequence $\lambda_{z,n}$, for which $|\lambda_{z,n}|\to 0$ when $n \to \infty$, satisfying  
\begin{equation*}
|f(z+\lambda_{z,n})-f(z)| \leq  e^{-\frac{C_1(2-s)^2}{p} \log^2 \left( \frac{1}{|\lambda_{z,n}|} \right)}
\end{equation*}
and 
\begin{equation*}
|\arg(f(z+\lambda_{z,n})-f(z))|\geq \frac{C_2(2-s)^2}{p} \log^{2} \left( \frac{1}{|\lambda_{z,n}|} \right).
\end{equation*}
\end{theorem}
The essential difference to positive results in Theorems \ref{VenytysT} and \ref{MSPEKTK} is the square in the term $2-s$, which corresponds to the difference of square root in the dimensional compression. Note that we can have both rotation and compression happening simultaneously, but the constants $C_1$ and $C_2$ depend on each other.\\
\\
Finally, we cover the question about expansion of dimension. That is, how big can a small set $A$ become under a homeomorphism with $p$-exponentially integrable distortion. This has previously been studied, in the form that we are interested in, by Clop and Herron in \cite{CH2} Corollary 3.1. They proved that if $f: \mathbb{C} \to \mathbb{C}$ is a homeomorphism with $p$-exponentially integrable distortion with $p>1$ and the gauge function is fixed as 
$$h(t)=\log^{\frac{-s(p-\epsilon)}{2}}\left( \frac{1}{t} \right),$$ 
then for each set $A$ such that $H^h(A)=0$ we have $H^s(f(A))=0$. This result follows from the upper bound  \eqref{PisteVen}, and hence it is sharp only when $s\to 0$. \\
\\
We can improve this result using the compression multifractal spectra for mappings with integrable distortion from \cite{H4}. Here the key point is that by a result of Gill \cite{GILL} we can couple the exponential integrability of the distortion of homeomorphism with the integrability of the distortion of the inverse.

\begin{theorem}\label{Korollaari}
Let  $f: \mathbb{C} \to \mathbb{C}$ be a homeomorphism with $p$-exponentially integrable distortion, where $p>1$, and fix $s\in (0,2)$ and the gauge function 
\begin{equation}\label{INTROGAUGE}
h(t)=\left( \frac{1}{\log\left( \frac{1}{t} \right)} \right)^{\frac{ps}{2-(s-\epsilon)}}.
\end{equation}
Under these assumptions if a set $A \subset \mathbb{C}$ satisfies $H^{h}(A)<\infty$ then $\dim(f(A))\leq s$. \\
\\
Moreover, this result is sharp in the sense that given any $\epsilon >0$ it is false for the gauge function 
\begin{equation*}
h(t)=\left( \frac{1}{\log\left( \frac{1}{t} \right)} \right)^{\frac{ps}{2-(s+\epsilon)}}.
\end{equation*}
The examples constructed to verify sharpness work for all $p>0$.
\end{theorem}
Note that when $s \to 0$ our gauge function \eqref{INTROGAUGE} converges towards the gauge function obtained by Clop and Herron, which is natural as Clop and Herron used the pointwise bound in their proof. \\
\\
Finally, let us say few words why the results in exponentially integrable case are not quite as sharp as when the distortion is assumed to be just integrable, see \cite{H4}. Our proofs are based on the modulus of path families, which is hard to calculate precisely and usually leaves a  constant coefficients in front of estimates. When the distortion is $p$-integrable these constants are not problematic as optimal bahaviour is studied, heuristically speaking, on exponential level. However, when the distortion is assumed to be of exponentially integrable order these constants directly affect the results, leading to what we believe to be non-optimal constants in \eqref{DIAMA}, \eqref{KIERTO} and \eqref{DIAMA2}. However, even with these complications we should get the right asymptotical behaviour when $s\to 2$, which makes the difference of square in term $2-s$ between the multifractal spectra and examples  perplexing.

\subsection*{Acknowledgements.} The author has been supported by the Academy of Finland project SA-1346562.

\section{Prerequisites} \label{se:prereq}
Let  $f:\Omega_1 \to \Omega_2$ be a sense-preserving homeomorphism between planar domains $\Omega_1, \Omega_2 \subset \mathbb{C}$. We say that $f$ is a \emph{$K$-quasiconformal mapping} for some $K \geq 1$ if $f \in W^{1,2}_{\text{loc}}$ and  the distortion inequality
\begin{equation*}
|Df(z)|^2 \leq KJ_f (z)
\end{equation*}
is satisfied almost everywhere. Here 
\begin{equation*}
|Df(z)|=\max \{|Df(z) e|: e\in \mathbb{C}, |e|=1\},
\end{equation*}
whereas $J_f(z)$ is the Jacobian of the mapping $f$ at the point $z$. 

Furthermore, we  say that $f$ has \emph{finite distortion} if the following conditions hold:
\begin{itemize}
	\item $f\in W_{\text{loc}}^{1,1}(\Omega_1)$
	\item $J_f(z)\in L^{1}_{\text{loc}}(\Omega_1)$
	\item $|Df(z)|^2\leq K(z)J_f(z) \qquad \text{almost everywhere in $\Omega_1$},$
\end{itemize}
for a measurable function $K(z)\geq 1$, which is finite almost everywhere. The smallest such function is denoted by $K_f(z)$ and called the distortion of $f$. Furthermore, we say that a mapping of finite distortion has \emph{$p$-exponentially integrable distortion} with a parameter $p>0$ if 
\begin{equation*}
e^{p K_{f}(z) } \in L^{1}_{\text{loc}.} 
\end{equation*}
This class of mappings generalizes quasiconformal mappings and has been intensively studied as many applications require that the distortion is allowed to blow up in a controlled manner. For a closer look on planar mappings of finite distortion see \cite{AIM}. In the present paper we only consider \textit{homeomorphic} mappings of finite distortion.
\medskip

Let $f: \mathbb{C} \to \mathbb{C}$ be a mapping of finite distortion and fix a point $z_0 \in \mathbb{C}$. In order to study the pointwise rotation of $f$ at the point $z_0$, we fix an argument $\theta\in[0,2\pi)$, and then look at how the quantity
$$\arg (f(z_0+te^{i\theta})-f(z_0))$$ 
changes as the parameter $t$ goes from 1 to a small $r$. This can also be understood as the winding of the path $f\left( [z_0+re^{i\theta}, z_0+e^{i\theta}] \right)$ around the point $f(z_0)$. As we are interested in the maximal pointwise spiraling we need to consider all directions $\theta$, leading to
\begin{equation}\label{spiraling}
\sup_{\theta \in [0,2\pi)} |\arg (f(z_0+re^{i\theta})-f(z_0)) - \arg (f(z_0+e^{i\theta})-f(z_0))  |.
\end{equation}

The maximal pointwise rotation is precisely the behavior of the above quantity \eqref{spiraling} when $r\to 0$. In this way, we say that the map $f$ \emph{spirals at the point $z_0$ with a rate function $g$}, where $g: [0, \infty) \to [0, \infty)$ is a decreasing continuous function, if 
\begin{equation}\label{Spiral rate}
\limsup_{r\to 0} \frac{\sup_{\theta \in [0,2\pi)} |\arg (f(z_0+re^{i\theta})-f(z_0)) - \arg (f(z_0+e^{i\theta})-f(z_0))  |}{g(r)} = C
\end{equation}
for some constant $C>0$. Finding maximal pointwise rotation for a given class of mappings equals finding the maximal spiraling rate for this class. Note that in \eqref{Spiral rate} we must use limit superior as the limit itself might not exist. Furthermore, for a given mapping $f$ there might be many different sequences $r_n \to 0$ along which it has profoundly different rotational behaviour. 
For more on these definitions we refer to  \cite{AIPS, H4}. 
\medskip 

In proofs of our theorems the \emph{modulus of path families} will play an important role. We provide here the main definitions, but interested reader can find a closer look at the topic in \cite{V}. The image of a line segment $I$ under a continuous mapping is called a path and a family of paths is denoted by $\Gamma$. Given a path family $\Gamma$ we say that a Borel measurable function $\rho:\mathbb{C} \to [0,\infty]$ is admissible with respect to  it if any locally rectifiable path $\gamma \in \Gamma$ satisfies
\begin{equation*}
\int_{\gamma} \rho(z) |dz| \geq 1.
\end{equation*}
Denote the modulus of the path family $\Gamma$ by $M(\Gamma)$ and define it by 
\begin{equation*}
M(\Gamma)= \inf_{\rho \text{ admissible}} \int_{\mathbb{C}} \rho^{2}(z) dA(z).
\end{equation*}
As an intuitive rule the modulus is large if the family $\Gamma$ has ``lots" of short paths, and is small if the paths are long and there is not ``many" of them. 

We will also need a \emph{weighted} version of the modulus. Let a weight function $\omega: \mathbb{C} \to [0, \infty)$, which in our case will always be the distortion function $K_f$, be  measurable and locally integrable.  Define the weighted modulus $M_{\omega}(\Gamma)$ by
\begin{equation*}
M_{\omega}(\Gamma)= \inf_{\rho \text{ admissible}} \int_{\mathbb{C}} \rho^{2}(z) \omega(z) dA(z).
\end{equation*}
One of the key ingrediens in our proofs is the modulus inequality
\begin{equation}\label{Modeq}
M(f(\Gamma)) \leq M_{K_f}(\Gamma)
\end{equation}
which holds for any mapping of finite distortion $f$ with locally integrable distortion, see \cite{H3}. \\
\\
Finally, we need the standard estimate for the modulus of path family $\Gamma$ connecting  two disjoint continua $E$ and $F$ in the form of Lemma 7.38 from \cite{V}, which gives the bound
\begin{equation}\label{MODLEMMA}
M(\Gamma) \geq  \tilde{C} \log \left( 1+ \frac{\min\{\text{diam}(E),\text{diam}(F)\}}{\text{dist} (E,F)} \right)
\end{equation}
where $\tilde{C}$ is some fixed constant that does not depend on sets $E$, $F$.

\section{Multifractal spectra}
In order to prove Theorems \ref{VenytysT} and \ref{MSPEKTK} we use ideas inspired by \cite{H4}, where similar results for mappings with $p$-integrable distortion were proved. The key idea is  to find small segments where mapping compresses or rotates strongly. These segments will be vital in finding good path families $\Gamma$ to use in the proofs.  Let us first consider compression.

\begin{lemma}\label{AUX0}
Fix parameters  $p>0$, $s \in (0,2)$ and $\epsilon >0$, and let $f$ be a homeomorphism. Assume that we can find a point $z\in \mathbb{C}$ with a sequence of complex numbers $\lambda_n$, such that the moduli $|\lambda_n|\to 0$ form a decreasing sequence, which satisfy
\begin{equation}\label{APU1}
|f(z+\lambda_n)-f(z)| \leq e^{- \frac{(2-s)\bar{C}}{2p} \log^2\left( \frac{1}{|\lambda_n|} \right)}
\end{equation}
for every $n$.  Then we can find a  sequence of complex numbers $\Lambda_n$, whose moduli satisfy $|\Lambda_n|=\frac{1}{e^{m_n}}$ for some increasing sequence of positive integers $m_n$, such that  from every line segment $$[z+\Lambda_n,z+e \Lambda_n]$$ 
we can find
points $x$ and $y$ satisfying 
\begin{equation}\label{AUX0y}
\frac{|f(x)-f(z)|}{|f(y)-f(z)|} \leq e^{-\frac{(2-\epsilon)(2-s)\bar{C}m_n}{2p}}.
\end{equation}
\end{lemma}
\proof: Without loss of generality fix $z=0$ and $f(z)=0$. Assume that we can't find line segments such that \eqref{AUX0y}  holds starting from some exponent $N \in \mathbb{N}$.  Denote $l=\min_{\theta\in [0,2\pi)} \left|f\left( \frac{e^{i\theta}}{e^{N}} \right)\right|$ and write $|f(\lambda_n)|$, where $|\lambda_n| \in \left[ \frac{1}{e^{M+1}}, \frac{1}{e^{M}} \right]$ with $M>N$, in the form
\begin{equation*}
|f(\lambda_n)|= \frac{|f(\lambda_n)|}{\left| f \left( \frac{\lambda_n}{|\lambda_n|} \frac{1}{e^{M}} \right) \right|} \cdot \frac{\left| f \left( \frac{\lambda_n}{|\lambda_n|} \frac{1}{e^{M}} \right) \right|}{\left| f \left( \frac{\lambda_n}{|\lambda_n|} \frac{1}{e^{M-1}} \right) \right|} \cdots \frac{\left| f \left( \frac{\lambda_n}{|\lambda_n|} \frac{1}{e^{N+1}} \right) \right|}{\left| f \left( \frac{\lambda_n}{|\lambda_n|} \frac{1}{e^{N}} \right) \right|} \cdot\left| f \left( \frac{\lambda_n}{|\lambda_n|} \frac{1}{e^{N}} \right) \right|.
\end{equation*}
Then for each division use the assumption that \eqref{AUX0y} does not hold and estimate the last term by $l$, which leads to
\begin{equation*}
|f(\lambda_n)| \geq l \exp\left( -\sum_{j=0}^{M+1}  \frac{(2-\epsilon)(2-s)\bar{C}j}{2p}\right)=l \exp\left(- \frac{(2-\epsilon)(2-s)\bar{C}(M+1)(M+2)}{4p}  \right).
\end{equation*}
On the other hand we can estimate using \eqref{APU1}
\begin{equation*}
|f(\lambda_n)| \leq e^{- \frac{(2-s)\bar{C}}{2p} \log^2\left( \frac{1}{|\lambda_n|} \right)} \leq  e^{- \frac{(2-s)\bar{C}(M-1)^2}{2p} 
}.
\end{equation*}
But now as $\epsilon>0$ is fixed  we get a contradiction when $M\to \infty$, that is, when $|\lambda_n| \to 0$. Hence the original claim holds and we can find such sequence $\Lambda_n$. \\
\\
Next we show a similar result for rotation.

\begin{figure}
\begin{center}
\includegraphics[scale=0.7]{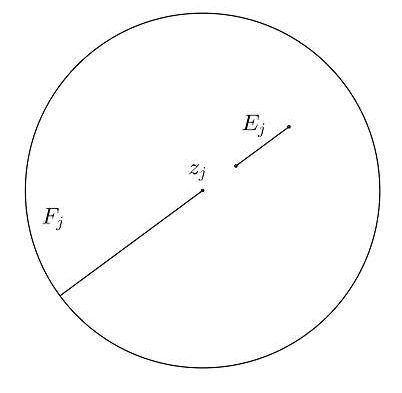}
\caption{The sets $E_j$ and $F_j$.}
\label{FIG1}
\end{center}
\end{figure}

\begin{lemma} \label{AUX1}
Fix parameters  $p>0$, $s \in (0,2)$ and $\epsilon>0$, and let $f$ be a homeomorphism. Assume that we can find a point $z\in \mathbb{C}$ with a  sequence of complex numbers $\lambda_n$, for which the moduli $|\lambda_n|\to 0$ form a decreasing sequence, and a branch of the argument which satisfy 
\begin{equation}\label{AUX1APU}
|\arg(f(z+\lambda_n)-f(z))|\geq  \frac{(2-s)\sqrt{\pi \tilde{C}}\bar{C} }{\sqrt{2} p} \log^2 \left( \frac{1}{|\lambda_n|} \right)
\end{equation}
for  every $n$. Then we can find a  sequence of complex numbers $\Lambda_n$, whose moduli satisfy $|\Lambda_n|=\frac{1}{e^{m_n}}$ for some increasing sequence of positive integers $m_n$, such that  from every line segment 
\begin{equation}\label{Lisäys}
[z+\Lambda_n,z+e \Lambda_n]
\end{equation}
  we can find
points $x$ and $y$ satisfying
\begin{equation}\label{AUX1y}
|\arg(f(z+x)-f(z))-\arg(f(z+y)-f(z))| \geq  \frac{(2-\epsilon)(2-s)\bar{C} \sqrt{\tilde{C} \pi}m_n}{\sqrt{2} p}.
\end{equation} 
\end{lemma}
Proof: We will again without loss of generality fix $z=0$ and $f(z)=0$. Assume that starting from some $N\in \mathbb{N}$ we can't find segment  \eqref{Lisäys} such that \eqref{AUX1y} is satisfied.  Denote $$\alpha= \max_{\theta\in [0,2\pi)} \left|\arg\left( f \left( \frac{e^{i\theta}}{e^N} \right) \right)\right|$$ and estimate $|\arg(f(\lambda_n))|$, where $|\lambda_n| \in \left[ \frac{1}{e^{M+1}}, \frac{1}{e^{M}} \right]$ and $M>N$, by 
\begin{equation*}
\begin{split}
|\arg(f(\lambda_n))| & \leq   \alpha+ \left|\arg\left( f\left(\frac{\lambda_n}{|\lambda_n|}\frac{1}{e^{N+1}}\right) \right)- \arg\left( f\left(\frac{\lambda_n}{|\lambda_n|}\frac{1}{e^{N}} \right) \right)\right| + \cdots \\
& + \left|\arg\left( f(\lambda_n) \right)- \arg\left( f \left(\frac{\lambda_n}{|\lambda_n|}\frac{1}{e^{M}} \right) \right) \right|.
\end{split}
\end{equation*}
Estimate each term using the assumption that \eqref{AUX1y} does not hold to obtain
\begin{equation*}
\begin{split}
|\arg(f(\lambda_n))| \leq \alpha + \frac{(2-\epsilon)(2-s)\bar{C} \sqrt{\tilde{C}\pi}}{\sqrt{2} p} \sum_{j=0}^{M+1}m_n= \alpha + \frac{(2-\epsilon)(2-s)\bar{C} \sqrt{\tilde{C}\pi}(M+1)(M+2)}{2 \cdot \sqrt{2} p}.
\end{split}
\end{equation*}
On the other hand we can use \eqref{AUX1APU} to estimate
\begin{equation*}
|\arg(f(\lambda_n))| \geq \frac{(2-s)\sqrt{\tilde{C}\pi}\bar{C}}{\sqrt{2} p} \log^2 \left( \frac{1}{|\lambda_n|} \right) \geq  \frac{(2-s)\sqrt{\tilde{C}\pi}\bar{C} (M-1)^2}{\sqrt{2} p}.
\end{equation*}
When $|\lambda_n|\to 0$ we see that $M\to \infty$, and thus the above estimates yield a contradiction. Hence the claim holds and we can always find such a sequence $\Lambda_n$. \\
\\
Lemmata \ref{AUX0} and \ref{AUX1} show that when studying compression or rotation we can focus on small segments without losing too much information. Our aim is to use these line segments as building blocks for path family, but first we must insulate them from each other. To be more precise, we want to surround these line segments with circle like sets as in figure \ref{FIG1}.  \\
\\
Our next task is to show that we can find many such building blocks that are disjoint from each other. To this end we use Lemma 3.4 from \cite{H4}, which ensures that we can find arbitrary small radii $$R_{n_0}=\frac{1}{e^{n_0}}$$ such that there exists many points $z_j\in A$ with the line segments of form $$\left[z_j+\Lambda_{z_j}, z_j+e\Lambda_{z_j} \right],$$ where $|\Lambda_{z_j}|=R_{n_0}$, that satisfy either the stretching condition \eqref{AUX0y} or the rotational condition  \eqref{AUX1y}. Moreover, we show that we can choose these points $z_j$  so that 
the circle like sets $F_j$, centered at the points $z_j$ with  radius of $e^2 R_{n_0}$, are disjoint and that each of them encircles the corresponding line segment $\left[z_j+\Lambda_{z_j}, z_j+e\Lambda_{z_j} \right]$.

\begin{lemma} \label{AUX2}
[Lemma 3.4 in \cite{H4}] Fix any $\delta \in (0,1)$, $s\in (0,2)$ and $\alpha>0$, and let $A \subset \mathbb{C}$. Assume that $\dim(A)=s$ and associate with every point $z\in A$ a decreasing sequence of radii $\{\delta^{k_{z,n}}\}_{n=1}^{\infty}$, where $k_{z,n}$ is a sequence of positive integers. Then for any given $s_0<s$ there exists radii $\delta^k$, which we can choose as small as we wish, such that we can find $\left\lfloor\left( \frac{1}{\alpha \delta^k} \right)^{s_0} \right\rfloor$ disjoint balls $\overline{B}(z_j,\alpha \delta^k)$, where $z_j \in A$ and $\delta^k \in \{\delta^{k_{z_j,n}}\}_{n=1}^{\infty}$ for every $j$.
\end{lemma}
Armed with these auxiliary lemmata, which are necessary for finding suitable path families $\Gamma$ to use in proofs, we are ready to move on to the main theorems.

\subsection{Proof of Theorem \ref{VenytysT}} Let us first concentrate on the compression multifractal spectra. Fix parameters  $p > 0$, $s\in (0,2)$ and, as we are interested in the Hausdorff dimension, we can assume $A \subset D$. Using Lemma \ref{AUX0} we can find for every point $z\in A$ a sequence of complex numbers $\Lambda_{z,n}$, whose moduli form a decreasing sequence 
\begin{equation} \label{ylimääräinenviite1}
\{|\Lambda_{z,n}|\}_{n=1}^{\infty}=\left\{\frac{1}{e^{m_{z,n}}} \right\}_{n=1}^{\infty},
\end{equation}
 such that the segments  $\left[z+ \Lambda_{z,n}, z+e\Lambda_{z,n} \right]$ satisfy the stretching condition \eqref{AUX0y}. Then using Lemma \ref{AUX2}, with the choices $\delta=\frac{1}{e}$ and $a=e^2$, we can find an arbitrary small radius $$R_{n_0}=\frac{1}{e^{n_0}}$$
for which there exists 
$$\left\lfloor\left( \frac{1}{e^{2} R_{n_0}} \right)^{\dim(A)-\bar{\epsilon}}\right\rfloor$$ 
disjoint circles with radius $e^2 R_{n_0}$ and centerpoints $z_j \in A$ that satisfy $R_{n_0} \in \{ |\Lambda_{z_{j},n}| \}_{n=1}^{\infty}$ for every $j$. Note that both the $\epsilon$ and the radius $R_{n_0}$ can be chosen arbitrary small. \\
\\
For any such centerpoint $z_j$ denote by $E_j$ the line segment  $\left[z_j+ \Lambda_{z_j},z_j+ e \Lambda_{z_j}  \right]$, where $\Lambda_{z_j}\in \{\Lambda_{z_j,n}\}_{n=1}^{\infty}$ and $| \Lambda_{z_j}|=R_{n_0}$, that satisfies the stretching condition \eqref{AUX0y}.  Such line segments clearly exists due to the choice of points $z_j$. Then denote 
\begin{equation*}
F_j=\left[z_j, z_j+ e^{i\pi}e^{2} \Lambda_{z_j} \right] \bigcup \partial B(z_j,e^{2}R_{n_0}).
\end{equation*}
For the illustration of the sets $E_j$ and $F_j$ see the figure \ref{FIG1}. Finally we define 
\begin{equation*}
E =\bigcup E_j \quad \text{and} \quad F=\bigcup F_j
\end{equation*}
and set $\Gamma$ to be the family of all paths connecting the sets $E$ and $F$. Note that, as each set $E_j$ is enclosed by the set $F_j$, from the modulus point of view we can think of our path family $\Gamma$ as the union of separate path families $\Gamma_j$ that connect $E_j$ to $F_j$ and consists of paths living inside the ball $B(z_j,e^2R_{n_0})$. \\
\\
As we want to use the modulus inequality \eqref{Modeq} we must  estimate the moduli $M_{K_f}(\Gamma)$ and $M(f(\Gamma))$. We will do this separately for each $\Gamma_j$ and start with the weighted modulus $M_{K_f}(\Gamma_j)$. To this end define a non-negative borel-measurable function
\begin{equation*}
\rho_{j}(z)=  \left\{
  \begin{array}{l l}
    \frac{1}{R_{n_0}} & \quad \text{if $\dist(z,E_j) < R_{n_0}$}\\
    0 & \quad \text{otherwise }
  \end{array} \right.
\end{equation*} 
which is clearly admissible for $\Gamma_j$. Hence we can estimate
\begin{equation}\label{Venytys1}
M_{K_f}(\Gamma_j) \leq \int_{\dist(z,E_j) \leq R_{n_0}} \rho_{0}^{2}(z) K_f(z)  dA(z)= \frac{1}{ R_{n_0}^2} \int_{\dist(z,E_j) < R_{n_0}}  K_f(z)  dA(z).
\end{equation}
On the other hand, to bound the modulus $M(f(\Gamma_j))$ from below we use the standard estimate in the form of \eqref{MODLEMMA} to obtain
\begin{equation*}
M(f(\Gamma_j )) \geq \tilde{C}\log\left(1 + \frac{\text{diam}(f(E_j))}{d(f(E_j), f(F_j))} \right).
\end{equation*}
We remind that the line segments $E_j$ satisfy the stretching condition \eqref{AUX0y} and hence
\begin{equation}\label{Venytys2}
M(f(\Gamma_j )) \geq \frac{(2-\epsilon)(2-s)\bar{C} \tilde{C} n_0}{2p}=\frac{(2-\epsilon)(2-s)\bar{C}\tilde{C} }{2p} \log\left( \frac{1}{R_{n_0}} \right).
\end{equation}
Combining \eqref{Venytys1} and \eqref{Venytys2} we get
\begin{equation}\label{Venytys3}
\frac{(2-\epsilon)(2-s)\bar{C} \tilde{C} }{(2(e-1)+\pi)2p} \log\left( \frac{1}{R_{n_0}} \right) \leq   \int_{\dist(z,E_j) < R_{n_0}} \frac{ K_f(z)}{ (2(e-1)+\pi)R_{n_0}^2}  dA(z).
\end{equation}
This  provides us a lower bound for  size of the distortion $K_f(z)$ inside the balls $B(z_j, e^2R_{n_0})$. To jump from this to the bound for the exponential of the distortion we use the Jensen's inequality with the convex function $e^{px}$ and the probability measure $\frac{dA(z)}{(2(e-1)+\pi)R_{n_0}^2}$ to obtain 
\begin{equation*}\label{Venytys4}
\frac{C}{R_{n_0}^2} \int_{\dist(z,E_j) < R_{n_0}} e^{pK_f(z)} dA(z) \geq \left(\frac{1}{R_{n_0}}\right)^{\frac{(2-\epsilon)(2-s)\bar{C}\tilde{C}}{2(2(e-1)+\pi)}}.
\end{equation*}
And since we have 
$$\left\lfloor\left( \frac{1}{e^{2} R_{n_0}} \right)^{\dim(A)-\bar{\epsilon}}\right\rfloor$$ 
such balls $B(z_j, e^2 R_{n_0})$ we see that in order for the distortion to be $p$-exponentially integrable the term 
\begin{equation*}
C R_{n_0}^2 \left(\frac{1}{R_{n_0}}\right)^{\frac{(2-\epsilon)(2-s)\bar{C}\tilde{C}}{2(2(e-1)+\pi)}} \left(\frac{1}{R_{n_0}}\right)^{\dim(A)-\bar{\epsilon}}
\end{equation*}
must stay bounded when $R_{n_0} \to 0$, which gives the desired bound for the dimension as we can choose $\epsilon$ and $\bar{\epsilon}$ as small as we want.

\subsection{Proof of Theorem \ref{MSPEKTK}} Let us then turn our attention to the rotational multifractal spectra. We again use the modulus inequality \eqref{Modeq} with the same sets $E$ and $F$ as in the compression case, now assuming that rotational condition \eqref{APU1} is satisfied for all segments $E_j$.  But unlike in stretching case we prove Theorem \ref{MSPEKTK} by contradiction, and thus we assume that there exists some fixed $\bar{\epsilon}>0$ such that  $$\dim(A)=2- \frac{\bar{C} \tilde{C} (2-s)}{\pi + 2(e-1)} + \bar{\epsilon}.$$
 For the weighted modulus of  path families $\Gamma_j$ we use the same estimate as in \eqref{Venytys1}, so we can turn our attention to the modulus of the image side $f(\Gamma_j)$. Here we use the  same method as in \cite{H2, H4}, using the fact that $f(E_j)$ and $f(F_j)$ must cycle around the point $f(z_j)$ at least $$n_{j,R_{n_0}}=\frac{|\arg(f(z_j+\lambda_{n_0})-f(z_j))|}{2\pi}-1$$ 
times. Going word by word through the details of estimating the modulus $M(f(\Gamma_j))$ in \cite{H4}, starting from (5.4) and finishing at equation (5.12), we obtain 
\begin{equation*}
M(f(\Gamma_j)) \geq 2\pi \frac{ n_{j,R_{n_0}}^{2}}{\log\left( \frac{\sup_{x\in E_j} |f(x)-f(z_j)|}{\inf_{x\in E_j} |f(x)-f(z_j)|} \right)}.
\end{equation*}
Then using the rotation bound \eqref{APU1} to estimate $n_{j,R_{n_0}}$ from below gives
\begin{equation}\label{Kierto1}
\begin{split}
M(f(\Gamma_j)) &\geq 2\pi  \left( \frac{(2-\epsilon)(2-s)\bar{C} \sqrt{\tilde{C}} n_{0}}{2^{ \frac{3}{2}} \sqrt{\pi} p} \right)^2 \frac{1}{\log\left( \frac{\sup_{x\in E_j} |f(x)-f(z_j)|}{\inf_{x\in E_j} |f(x)-f(z_j)|} \right)}\\
& = \frac{\left( (2-\epsilon)(2-s)\bar{C} \sqrt{\tilde{C}} n_{0}\right)^2}{2^2 p^2} \frac{1}{\log\left( \frac{\sup_{x\in E_j} |f(x)-f(z_j)|}{\inf_{x\in E_j} |f(x)-f(z_j)|} \right)} .
\end{split}
\end{equation}
Hence, combining estimates \eqref{Venytys1} and \eqref{Kierto1} with \eqref{Modeq} we get
\begin{equation*}
\frac{1}{2(e-1)+\pi}\frac{\left( (2-\epsilon)(2-s)\bar{C} \sqrt{\tilde{C}} n_{0}\right)^2}{2^2 p^2} \frac{1}{\log\left( \frac{\sup_{x\in E_j} |f(x)-f(z_j)|}{\inf_{x\in E_j} |f(x)-f(z_j)|} \right)} \leq   \int_{\dist(z,E_j) < R_{n_0}} \frac{ K_f(z)}{ (2(e-1)+\pi)R_{n_0}^2}  dA(z).
\end{equation*}
Here we use Jensen's inequality with the convex function $e^{px}$ and the probability measure $\frac{dA(z)}{(2(e-1)+\pi)R_{n_0}^2}$ together with the above estimate to obtain 
\begin{equation*}
\frac{C}{R_{n_0}^2} \int_{\dist(z,E_j) < R_{n_0}} e^{pK_f(z)} dA(z) \geq \exp\left( \frac{\left( (2-\epsilon)(2-s)\bar{C} \sqrt{\tilde{C}} n_{0}\right)^2}{2^2 p (2(e-1)+\pi)} \frac{1}{\log\left( \frac{\sup_{x\in E_j} |f(x)-f(z_j)|}{\inf_{x\in E_j} |f(x)-f(z_j)|} \right)} \right),
\end{equation*}
which holds for an arbitrary $j$. Hence, summing over all $j$ we arrive at 
\begin{equation}\label{Kierto2}
\int_{\mathbb{D}} e^{p K_f(z)} d A(z) \geq C R_{n_0}^2 \sum_{j} \exp\left( \frac{\left( (2-\epsilon)(2-s)\bar{C} \sqrt{\tilde{C}} n_{0}\right)^2}{2^2 p (2(e-1)+\pi)} \frac{1}{\log\left( \frac{\sup_{x\in E_j} |f(x)-f(z_j)|}{\inf_{x\in E_j} |f(x)-f(z_j)|} \right)} \right).
\end{equation}
Here the terms $$\frac{1}{\log\left( \frac{\sup_{x\in E_j} |f(x)-f(z_j)|}{\inf_{x\in E_j} |f(x)-f(z_j)|} \right)}$$ inside the sum make it hard to estimate \eqref{Kierto2} directly. Instead we prove the following auxiliary result, which shows that in most cases these terms won't be too big.
\begin{lemma}\label{AUX3}
Under the assumptions made in Theorem \ref{MSPEKTK} and during its proof, especially that for some fixed $\bar{\epsilon}>0$ we have  $$\dim(A)=2- \frac{\bar{C} \tilde{C} (2-s)}{\pi + 2(e-1)} + \bar{\epsilon},$$
there is more than  $$\frac{1}{2} \left\lfloor\left( \frac{1}{e^{2} R_{n_0}} \right)^{\dim(A)-\epsilon}\right\rfloor$$ 
segments $E_j$ that satisfy 
\begin{equation}\label{TAPU4}
\frac{\sup_{x\in E_j} |f(x)-f(z_j)|}{\inf_{x\in E_j} |f(x)-f(z_j)|} \leq e^{\frac{\bar{C}(2-s)n_{0}}{p}}
\end{equation}
whenever $n_{0} $ is big enough.
\end{lemma}
Proof of Lemma \ref{AUX3}: We will prove the lemma using contradiction, so assume that we find at least 
\begin{equation}\label{TAPU1}
\frac{1}{2} \left\lfloor\left( \frac{1}{e^{2} R_{n_0}} \right)^{\dim(A)-\epsilon}\right\rfloor
\end{equation}
segments $E_j$ such that 
\begin{equation}\label{TAPU2}
\frac{\sup_{x\in E_j} |f(x)-f(z_j)|}{\inf_{x\in E_j} |f(x)-f(z_j)|} \geq e^{\frac{\bar{C}(2-s)n_{0}}{p}}.
\end{equation}
Let us restrict families $E$ and $F$ to these sets and proceed as in the proof of the compression multifractal spectra. That is,  we again use the modulus inequality \eqref{Modeq} for each $\Gamma_j$ and estimate as in \eqref{Venytys1} that 
\begin{equation*}
M_{K_f}(\Gamma_j)  \leq  \frac{1}{ R_{n_0}^2} \int_{\dist(z,E_j) < R_{n_0}}  K_f(z)  dA(z)
\end{equation*}
and use the standard modulus estimate \eqref{MODLEMMA} together with the assumption \eqref{TAPU2} leading to
\begin{equation*}
M(f(\Gamma_j)) \geq \tilde{C} \log\left( \frac{\sup_{x\in E_j} |f(x)-f(z_j)|}{\inf_{x\in E_j} |f(x)-f(z_j)|}   \right) \geq \frac{\bar{C}\tilde{C} (2-s)}{p} \log \left( \frac{1}{R_{n_0}} \right).
\end{equation*}
The modulus inequality \eqref{Modeq} used with the above estimates and the Jensen's inequality, where we again choose the convex function $e^{px}$ and the probability measure $\frac{dA(z)}{(2(e-1)+\pi)R_{n_0}^2}$, gives us an estimate
\begin{equation*}
\frac{C}{R_{n_0}^2} \int_{\dist(z,E_j) < R_{n_0}}  e^{pK_f(z)} dA(z) \geq  \left( \frac{1}{R_{n_0}} \right)^{\frac{\bar{C}\tilde{C} (2-s)}{\pi+2(e-1)}}.
\end{equation*}
And finally going over all selected $j$, see \eqref{TAPU1}, we have 
\begin{equation}\label{TAPU3}
\int_{\mathbb{D}} e^{pK_f(z)} dA(z) \geq C R_{n_0}^2 \left( \frac{1}{R_{n_0}} \right)^{\frac{\bar{C} \tilde{C} (2-s)}{\pi+2(e-1)}}  \left( \frac{1}{ R_{n_0}} \right)^{\dim(A)-\epsilon}
\end{equation}
where the right hand side is, by the assumption on the dimension of $A$ , of order $$ \left( \frac{1}{ R_{n_0}} \right)^{\bar{\epsilon}-\epsilon},$$
where $\bar{\epsilon}>0$ was fixed and $\epsilon>0$ can be made arbitrary small. Thus we note that since the left hand side of \eqref{TAPU3} is bounded from above for an arbitrary fixed $f$ we get a contradiction when $n_0$ is big. Thus the claim of finding at least \eqref{TAPU1} segments was false and the Lemma \ref{AUX3} has been proven. \\
\\
We can now turn our attention back to the proof of the rotational multifractal spectra. Let us from now on consider only  segments $E_j$ that satisfy \eqref{TAPU4}. From Lemma \ref{AUX3} we know, since we can assume $n_0$ to be as big as we want, that there is at least $$\frac{1}{2} \left\lfloor\left( \frac{1}{e^{2} R_{n_0}} \right)^{\dim(A)-\epsilon}\right\rfloor$$ such segments. Thus we can continue the estimate \eqref{Kierto2} by 
\begin{equation*}
\begin{split}
\int_{\mathbb{D}} e^{p K_f(z)} d A(z) &\geq C R_{n_0}^2  \sum_{j}  \exp\left( \frac{\left( (2-\epsilon)(2-s)\bar{C} \sqrt{\tilde{C}} n_{0}\right)^2}{2^2 p (2(e-1)+\pi)}  \frac{p}{\bar{C} (2-s)n_{0}} \right) \\
& =C R_{n_0}^2  \sum_{j}  \exp\left( \frac{(1-\tilde{\epsilon})\bar{C} \tilde{C} (2-s)n_0} {2(e-1)+\pi} \right) \\
&=C R_{n_0}^2  \sum_{j}  \exp\left(  \frac{(1-\tilde{\epsilon})\bar{C} \tilde{C} (2-s)} {2(e-1)+\pi} \log \left( \frac{1}{R_{n_0}}\right)\right) \\
& \geq  C R_{n_0}^2   \left( \frac{1}{R_{n_0}}\right)^{\frac{(1-\tilde{\epsilon})\bar{C} \tilde{C} (2-s)} {2(e-1)+\pi}} \left( \frac{1}{R_{n_0}}\right)^{\dim(A)-\epsilon},
\end{split}
\end{equation*}
where $\tilde{\epsilon} \to 0$ when $\epsilon \to 0$. But the above estimate shows that if 
$$\dim(A)=2- \frac{\bar{C} \tilde{C} (2-s)}{\pi + 2(e-1)} + \bar{\epsilon},$$
for fixed $\bar{\epsilon}>0$, we get a contradiction with $f$ having $p$-exponentially integrable distortion as we can  choose $\epsilon$ as small as we want and $n_0$ as big as we want. This finishes the proof of Theorem \ref{MSPEKTK}.

\section{Dimension compression and expansion}
Let us first concentrate on  compression of dimension and use the compression multifractal spectra to prove Theorem \ref{DimensioT}. By substituting  
$$s= 2- \frac{2(e-1)+\pi}{\bar{C}\tilde{C}}(2-\bar{s})$$
in Theorem \ref{VenytysT} we obtain the same result with the assumption \eqref{MSVEN} on compression  replaced by 
\begin{equation}\label{VVenytys}
|f(z+\lambda_n)-f(z)| \leq e^{-\frac{(2(e-1)+\pi)(2-\bar{s})}{2p\tilde{C}} \log^2\left( \frac{1}{|\lambda_n|} \right)}
\end{equation}
and the bound for the dimension in form
\begin{equation}\label{VVVenytys}
\dim(A) \leq \bar{s}.
\end{equation}
Here we immediately note that as $s\in (0,2)$ we have $$\bar{s} \in \left( 2- \frac{2 \bar{C} \tilde{C}}{2(e-1)+\pi}  ,2\right),$$
which gives the range for the dimension in which our result improves that of Zapadinskaya \cite{Zap}.  \\
\\
 Fix $\epsilon>0$ and the gauge function
\begin{equation}\label{Gauge}
h(t)= e^{-\log^{\frac{1}{2}} \left( \left(\frac{1}{t}\right)^{ \frac{2p \tilde{C} \bar{s}^2}{(2(e-1)+\pi)(2-\bar{s}+ \epsilon)} } \right)},
\end{equation}
and let $f: \mathbb{C} \to \mathbb{C}$ be a homeomorphism with $p$-exponentially integrable distortion. Our aim is to show that each $A\subset \mathbb{C}$ for which $H^{h}\left( f(A) \right)=0$ satisfies $H^{\bar{s}}(A)=0$. \\
\\
We can assume $A\subset \mathbb{D}$ without loss of regularity as countable union of sets with measure zero has measure zero. Denote by $A_1 \subset A$ the set of points for which there exists a sequence $\lambda_{z,n}$ such that $|\lambda_{z,n}|\to 0$ when $n\to \infty$ and
\begin{equation}\label{Vieläviite}
|f(z+\lambda_{z,n})-f(z)| \leq  e^{-\frac{(2(e-1)+\pi)(2-\bar{s}+ \epsilon)}{2p\tilde{C}} \log^2\left( \frac{1}{|\lambda_{z,n}|} \right)}.
\end{equation}
From Theorem \ref{VenytysT}, and the remarks \eqref{VVenytys}, \eqref{VVVenytys}, we know that $\dim(A_1) < \bar{s}$ and thus $H^{\bar{s}}(A_1)=0$. Hence we can concentrate on the set $A_2= A \setminus A_1$. \\
\\
Next we note that since for points $z\in A_2$ there does not exist a sequence $\lambda_{z,n}$ such that \eqref{Vieläviite} is satisfied we know that there exists a radius $r_z>0$ such that 
\begin{equation*}
|f(z+h)-f(z)| \geq  e^{-\frac{(2(e-1)+\pi)(2-\bar{s}+ \epsilon)}{2p\tilde{C}} \log^2\left( \frac{1}{|h|} \right)}
\end{equation*}
whenever $|h| < r_z$. Hence we get that for each $x=f(z)$, when $z\in A_2$, there exists radius $r_x >0$ such that 
\begin{equation}\label{KVenytys}
 |f^{-1}(x+h)-f^{-1}(x)|  \leq e^{-\log^{\frac{1}{2}} \left( \left(\frac{1}{|h|}\right)^{ \frac{2p \tilde{C} }{(2(e-1)+\pi)(2-\bar{s}+ \epsilon)} } \right)}
\end{equation}
whenever $|h|<r_x$. Denote 
$$A_{2,k}= \left\{ z\in A_2: r_x \geq \frac{1}{k} \qquad \text{where $x=f(z)$}\right\}$$
and note that $A_2 = \cup_{k=1}^{\infty} A_{2,k}$. Hence it is enough to show that $H^{\bar{s}}(A_{2,k})=0$ for an arbitrary $k\in \mathbb{N}$. \\
\\
Fix $k\in \mathbb{N}$ and note that since we assumed that $H^{h}(f(A))=0$ also $H^{h}(f(A_{2,k}))=0$. Hence we can cover the set $f(A_{2,k})$ with balls $B(x_i,r_{x,i})$, where $r_{x,i}<\frac{1}{k}$, $x_i \in f(A_2)$ and 
\begin{equation}\label{TAASYKSI}
\sum_{i} h(r_{x_i,i}) \leq \epsilon_1
\end{equation}
for an arbitrary predefined $\epsilon_1$. Next we estimate the Hausdorff measure of the set $A_{2,k}$ by covering it using sets $f^{-1}\left( B(x_i,r_{x,i}) \right)$, whose diameter can be controlled by \eqref{KVenytys}. To this end we calculate
\begin{equation*}
\begin{split}
\sum_{i} \left(  \dim \left( f^{-1}\left( B(x_i,r_{x_i,i}) \right) \right)  \right)^{\bar{s}} & \leq 2^{\bar{s}} \sum_{i}  e^{-\bar{s}\log^{\frac{1}{2}} \left( \left(\frac{1}{r_{x_i,i}}\right)^{ \frac{2p \tilde{C} }{(2(e-1)+\pi)(2-\bar{s}+ \epsilon)} } \right)} \\
& \leq 2^{\bar{s}} \sum_{i} h(r_{x_i,i}) \leq 4\epsilon_1,
\end{split}
\end{equation*}
where the last inequality follows from \eqref{TAASYKSI}. As we can choose $\epsilon_1$ and the radii $r_{x_i,i}$  arbitrary small this proves that $H^{\bar{s}} (A_{2,k})=0$ for any fixed $k \in \mathbb{N}$, and thus also $H^{\bar{s}}(A_2)=0$ finishing the proof. \\
\\
Theorem \ref{DimensioT} shows that the gauge function \eqref{Gauge} measuring contraction is of different order when $\bar{s}\to 2$  than previously obtained result by Zapadinskaya in \cite{Zap}. To the other direction examples by Zapadinskaya and Clop and Herron show that with the term $(2-\bar{s})^2$ we would get a sharp result. It is an interesting question if the term $2-\bar{s}$ or its square is the optimal one?  
\subsection{Proof of Theorem \ref{Korollaari}} To finish this section we consider the expansion of dimension for mappings with $p$-exponentially integrable distortion in the form of Theorem \ref{Korollaari}.\\
\\
Let $f:\mathbb{C} \to \mathbb{C}$ be a homeomorphism with $p$-exponentially integrable distortion. Then by the result of Gill \cite{GILL} $f^{-1}$ is a homeomorphism with $q$-integrable distortion for any $q<p$, which lets us couple expansion of homeomorphisms with exponentially integrable distortion and contraction of homeomorphisms with integrable distortion. \\
\\
Let us assume, in contradiction with the Theorem \ref{Korollaari}, that we can find parameters $\epsilon, \tilde{\epsilon}>0$, $p>1$ a set $E\subset \mathbb{C}$ and a homeomorphism $f:\mathbb{C} \to \mathbb{C}$ with $p$-exponentially integrable distortion such that $H^{h}(E)< \infty$, where $$h(t)= \left( \frac{1}{\log\left( \frac{1}{t} \right)} \right)^{\frac{
ps}{2-(s-\epsilon)}},$$ and $\diam(f(E))=s+\tilde{\epsilon}$. \\
\\
Since the inverse $f^{-1}$  has $(p-\bar{\epsilon})$-integrable distortion, for an arbitrary small $\bar{\epsilon}>0$, we can use Theorem 1.4 from \cite{H4} with the choice of $\epsilon=\frac{\tilde{\epsilon}}{4}$. From the above assumptions it follows that $H^{s+\frac{\tilde{\epsilon}}{2}}(f(E))>0$, and hence Theorem 1.4 used to the set $f(E)$  gives  $H^{\bar{h}}(E)>0$  with the gauge function $$\bar{h}(t)=\left( \frac{1}{\log\left( \frac{1}{t} \right)} \right)^{\frac{(p-\bar{\epsilon})\left(s+\frac{\tilde{\epsilon}}{2}\right)}{2-\left(s+ \frac{\tilde{\epsilon}}{4}\right)}}.$$
However, the set $E$ and the gauge functions $h$ and $\bar{h}$ must  satisfy  the basic inequality
\begin{equation*}
H^{\bar{h}}(E) \leq \limsup_{t\to 0} \frac{\bar{h}(t)}{h(t)} H^{h}(E).
\end{equation*} 
But now choosing $\bar{\epsilon}$  small enough leads to a contradiction with the original assumption $H^h(E)< \infty$ as $$\frac{\bar{h}(t)}{h(t)} \to 0$$ when $t\to 0$, which finishes the proof of the first part of the Theorem \ref{Korollaari}. The sharpness will be proven in the next chapter along with other examples. Note that we don't get sharpness straight from the examples proving sharpness of Theorem 1.4 in \cite{H4} as the inverses of those mappings don't have the correct regularity for the distortion.

\section{Examples}
\begin{figure}
\begin{center}
\includegraphics[scale=1]{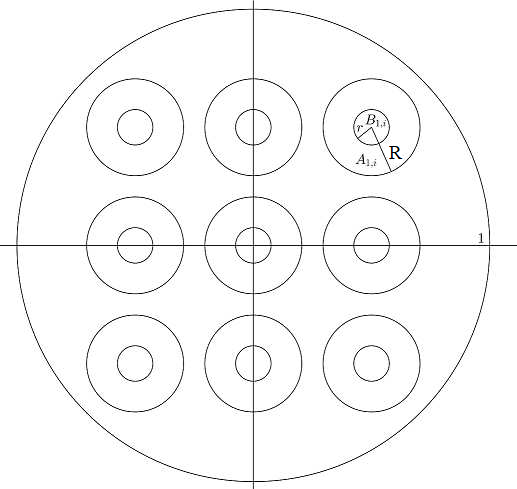}
\caption{The first level Balls $B_{1,i}$ and annuli $A_{1,i}$. }
\label{FIG2}
\end{center}
\end{figure}

We first  provide examples towards sharpness of the multifractal spectras. These are similar to classical examples using families of nested annuli forming cantor set at the limit, which in absense of rotation closely resemble examples presented, among others, by Zapadinskaya in \cite{Zap} or Clop and Herron, example 4.2 in \cite{CH}. For details on construction with rotation one can also see similar example in \cite{H4}. We will be somewhat careless with constants in the example as the main interest lies in asymptotical behaviour when $s\to 2$. \\
\\
At the first level of our construction choose $M$ balls $B_{1,i}$ with the radius $$r=\frac{1}{M^{s}}.$$ 
We pack the balls $B_{1,i}$ inside the unit disk in a square grid so that their distance from each other and from the boundary of the unit disk is at least $2er$, see figure \ref{FIG2}. In order for this packing to be possible we need to choose $M$ big enough. To be more precise we need to make sure that $M^{\frac{2-s}{2s}}<\sqrt{2}$ so that we can fit our square grid inside the unit disk, which gives the bound $M\geq C^{\frac{2s}{2-s}}$.\\
\\

\noindent Then we construct annuli out of the balls by fixing $A_{1,i}=eB_{1,i} \setminus B_{1,i}$ and note they are disjoint and stay inside the unit disk. We denote the family of first level ball by $\mathbb{B}_1$ and the family of annuli by $\mathbb{A}_1$. \\
\\
Given a similarity mapping $\phi_i : \mathbb{D} \to B_{1,i}$, which we will soon define, we create the family of the second level balls by
\begin{equation*}
\mathbb{B}_2=\bigcup_{i=1}^{M} \phi_{i}(\mathbb{B}_1),
\end{equation*}
and turn them to annuli with the family
\begin{equation*}
\mathbb{A}_2=\bigcup_{i=1}^{M} \phi_{i}(\mathbb{A}_1).
\end{equation*}
Note that the second level family of annuli has the inner radii of $r^2$ and the outer radii of $er^2$. \\
\\
We continue in iterative manner and define 
\begin{equation}
A=\bigcap_{n=1}^{\infty} \bigcup_{B\in \mathbb{B}_n} \overline{B},
\end{equation} 
which is a self similar Cantor set of the plane. The Hausdorff dimension of the set $A$ can be computed to be $s$ by checking that it satisfies the equation $$M r^s=1.$$ 
Then we define for an arbitrary annulus $B \setminus  \frac{1}{e} B$, with $B=B(a,R)$, and an arbitrary $K\geq 1$  the radial map
\begin{equation}\label{rakennuspala}
\psi_{B,K}(z)=  \left\{
  \begin{array}{l l}
    z & \quad \text{ if $z\notin B(a,R)$}\\
   a+ R\frac{z-a}{|z-a|}\left| \frac{z-a}{R}\right|^{(c_1+ic_2)\bar{K}} & \quad  \text{ if $z\in B(a,R) \setminus B\left(a,\frac{R}{e}\right)$} \\
   a+e^{-c_1\bar{K}+1}e^{-ic_2\bar{K}}(z-a) & \quad \text{ if $z\in B\left(a,\frac{R}{e}\right)$}
  \end{array} \right.
\end{equation} 
where $c_1>0 $ and $c_2 \geq 0$, which will serve as the building blocks for the mapping $f$. Note that the map $\psi_{B,K}(z)$ is $C_{c_1,c_2}\bar{K}$ quasiconformal and similarity, and thus conformal, outside the annulus $B \setminus  \frac{1}{e} B$.   Let us then fix  $\epsilon>0$ and set $$\bar{K}_n=\frac{1}{pC_{c_1,c_2}} \log \left( M^{\frac{(2-s-\epsilon)n}{s}} \right)$$ so that the distortion at the level $n$ gets the form
\begin{equation}\label{KMÄÄR}
K_{n}= \frac{1}{p} \log \left( M^{\frac{(2-s-\epsilon)n}{s}} \right),
\end{equation}
which we will use as the distortion at the corresponding step of the construction. \\
\\
Let us then begin the construction of our mapping $f$.  As the first step we define 
\begin{equation*}
f_{1}(z)= \left\{
  \begin{array}{l l}
    z & \quad \text{ if $z\notin \bigcup_{B\in \mathbb{A}_1} B$}\\
   \psi_{B,K_1}(z) & \quad \text{ if $z\in B$ with $B \in \mathbb{A}_1$}
  \end{array} \right.
\end{equation*} 
Then we assume that the previous step $f_{n-1}$ has been defined and set
\begin{equation*}
f_{n}(z)= \left\{
  \begin{array}{l l}
    f_{n-1}(z) & \quad \text{ if $z\notin \bigcup_{B\in \mathbb{A}_n} B$}\\
   \psi_{f_{n-1}(B),K_n}(f_{n-1}(z)) & \quad \text{ if $z\in B$ with $B \in \mathbb{A}_n$}
  \end{array} \right.
\end{equation*} 
Since $f_{n-1}$ is conformal inside the family $\mathbb{A}_n$ we see that $f_n$ is $K_{i}$-quasiconformal inside the the level $i$ annuli, where $i\in \{1,2,...,n\}$, and conformal elsewhere. 
 \\
\\
The sequence $f_n$ form a cauchy sequence and thus there exists a limit map $$f=\lim_{n\to \infty} f_n$$ which is clearly a homeomorphism. Moreover, the limit is differentiable almost everywhere as it is differentiable outside of the set $A$ and boundaries of annuli $A_{n,i}$. Classical calculations for the radial mappings $\psi_{B,K}$ give that $$|Df(z)| \leq C_{c_1,c_2}K_n$$
when $z$ is inside an annuli from level $n$ and that $|Df(z)|<1$ if $f$ is differentiable at $z$ which is not inside any annuli. Hence in order to show that $Df \in L^{1}_{\text{loc}}(\mathbb{C})$ it is enough to estimate the size of \eqref{KMÄÄR} as follows
\begin{equation}\label{KPINTLASKU}
\begin{split}
\sum_{n=1}^{\infty} |\mathbb{A}_n|e^{pK_n} & \leq c\sum_{n=1}^{\infty} M^n r^{2n}e^{pK_n} \\
& = c \sum_{n=1}^{\infty} \frac{M^n}{M^{\frac{2n}{s}}}M^{\frac{n(2-s-\epsilon)}{s}} \\
& = c \sum_{n=1}^{\infty} \frac{1}{M^{ \frac{\epsilon}{s}n}} < \infty
\end{split}
\end{equation}
 since local integrability follows from exponential integrability for any fixed $p$. \\
\\
As we arranged the balls $B_{1,i}$  in a square grid, see figure \ref{FIG2},  it follows that for almost every line parallel to the coordinate axis we can find a neighbourhood where the mapping $f$ coincides with the quasiconformal mappings $f_n$, for every big enough $n$. Thus the mapping $f$ is absolutely continuous for almost every line parallel to the coordinate axes. This together with  $Df(z) \in L^{1}_{\text{loc}}(\mathbb{C})$ implies that $f\in W^{1,1}_{\text{loc}}(\mathbb{C})$. Furthermore, since $f$ is a homeomorphism, this also shows that $J_f(z)\in L^{1}_{\text{loc}}(\mathbb{C})$. Thus the mapping $f$ is a mapping of finite distortion, with the distortion 
\begin{equation*}
K_{f}(z)= \left\{
  \begin{array}{l l}
    K_{n} & \quad \text{ if $z\in A_{n,i}$ for some $i\in \{1,2,...,M^n\}$}\\
    1 & \quad  \text{ otherwise} 
  \end{array} \right.
\end{equation*}
that is $p$-exponentially integrable due to the estimate \eqref{KPINTLASKU}. \\
\\
Hence we can move on to check rotation and compression at points $z\in A$. \newline \newline 
For rotation one can check in a similar manner as in \cite{H4} that it essentially, apart from small error term, comes from crossing the annuli $A_n$. To be more precise, one can check using methods in \cite{H1} and \cite{H4}  that for any $z\in A$ we can find sequence of complex numbers  $\lambda_{z,n}$ such that $|\lambda_{z,n}|=r^n$ and 
\begin{equation*}
\begin{split}
|\arg(f(z+\lambda_{z,n})-f(z))| & \geq \sum_{j=1}^{n} \arg\left( e^{ic_2 \bar{K}_{j}} \right)  -5\pi n \\
& \geq \sum_{j=1}^{n} \frac{c_2}{C_{c_1,c_2}p} \log \left( M^{\frac{(2-s-\epsilon)j}{s}} \right)  -5\pi n   \\
& \geq \frac{c_2 (2-s-\epsilon)\log(M)}{C_{c_1,c_2}p}\sum_{j=1}^{n}  j -5\pi n   \\
& \geq \frac{c_2 (2-s-\epsilon)^2 n(n+1)}{2C_{c_1,c_2}p}-5\pi n
\end{split}
\end{equation*}
where in the last inequality we used the bound for $M$. When $n\to \infty$ we get 
\begin{equation*}
|\arg(f(z+\lambda_{z,n})-f(z))|  \geq \frac{c_2 (2-s-\epsilon)^2 n^2}{2C_{c_1,c_2}p}=\frac{C(2-s)^2}{p} \log^{2} \left( \frac{1}{|\lambda_{z,n}|} \right).
\end{equation*}
For compression we can find for any $z\in A$ the  unique sequence of nested balls $B_{n,i_n}$ such that $$z= \bigcap_{n=1}^{\infty} B_{n,i_n}.$$ Since the point $z$ lies inside the balls $ B_{n,i_n}$ we can find at every step $n$ a complex number $\lambda_{z,n}$  that satisfies  $z+\lambda_{z,n} \in \partial B_{n,i_n} $  and  $$|\lambda_{z,n}|=r\left( B_{n,i_n} \right)=r^n.$$
For this sequence $\lambda_{z,n}$ we get 
\begin{equation*}
|f(z+\lambda_{z,n})-f(z)|=e^{-c_1(\bar{K_1}+ \cdots +\bar{K_n})+n}|\lambda_{z,n}|
\end{equation*}
which simplifies, for big $n$, to 
\begin{equation*}
|f(z+\lambda_{z,n})-f(z)|\leq e^{\frac{-c_1 (2-s-\epsilon)^2 n^2}{2C_{c_1,c_2}p}}=e^{-\frac{C(2-s)^2}{p} \log^2 \left( \frac{1}{|\lambda_{z,n}|} \right)}
\end{equation*}
and hence we get the compression and rotation claimed by Theorem \ref{Esimerkki}. As mentioned before the main difference between Theorems \ref{VenytysT} and \ref{MSPEKTK} is that instead of $2-s$ we have $(2-s)^2$ as the leading term when $s\to 2$. 
\subsection{Proof of Theorem \ref{Korollaari}}
To finish we show the sharpness of the Theorem \ref{Korollaari}. Hence we want to find a mapping $f: \mathbb{C} \to \mathbb{C}$ with $p$-exponentially integrable distortion and set $E$ which is small in the sense of the gauge function $$h(t)=\frac{1}{\log^{\frac{ps}{2-s-\bar{\epsilon}}}\left( \frac{1}{t} \right)},$$ but whose image under the mapping $f$ is $s-\epsilon$ dimensional, where $\epsilon $ and $\bar{\epsilon}$ are independent of each other and can be chosen to be arbitrary small. \\
\\
The basic idea will be similar as in the previous construction, but we must be more carefull in estimating distortion. Thus instead of the building block \eqref{rakennuspala} we fix $p>0$ and define for any given annuli $A=B(a,\tilde{R}) \setminus B(a,\tilde{r})$ and $\epsilon >0$ the radial mapping 
\begin{equation}\label{rakennuspala2}
\psi_{A}(z)=  \left\{
  \begin{array}{l l}
    z & \quad \text{ if $z\notin B(a,\tilde{R})$}\\
   a+ \frac{z-a}{|z-a|} \frac{\tilde{R} \log^{\frac{p}{2-\epsilon}} \left( \frac{1}{\tilde{R}} \right)}{\log^{\frac{p}{2-\epsilon}} \left( \frac{1}{|z-a|} \right)} & \quad  \text{ if $z\in B(a,\tilde{R}) \setminus B\left(a,\tilde{r} \right)$} \\
   a+ \frac{\tilde{R} \log^{\frac{p}{2-\epsilon}} \left( \frac{1}{\tilde{R}} \right)}{\tilde{r} \log^{\frac{p}{2-\epsilon}} \left( \frac{1}{\tilde{r}} \right)} (z-a) & \quad \text{ if $z\in B\left(a,\tilde{r}\right).$}
  \end{array} \right.
\end{equation} 
Note that this building block is still a similarity, and thus conformal, outside of the annulus $A$ and  quasiconformal inside it. But now the value of the distortion depends on $|z|$ and can be calculated to be
\begin{equation}\label{KKKK}
 K_{\psi_{A}}(z)= \frac{2-\epsilon}{p} \log \left( \frac{1}{|z-a|} \right), 
\end{equation}
when $z\in A$, for details see \cite{CaHe}. For an annuli $A$ we denote by $B_o$ its outer ball $B(a,\tilde{R})$ and by $B_i$ its inner ball $B(a,\tilde{r})$. \\
\\
Let us then fix $s\in (0,2)$,  a small number $r$, a bigger number $R=Cr^{\frac{s}{2}}$ and $M=\frac{1}{r^s}$, where we choose $r$ and constant $C$ so that $M$ is an integer and  we can fit $M$ balls with radius $R$ inside the unit disc. Finally fix the radii that we will use in our building blocks as
\begin{equation*}
\bar{r}_n = e^{- \left( \frac{1}{r} \right)^{\frac{(2-s)n}{p}}} \quad \text{and} \quad \bar{R}_n=R\bar{r}_{n-1}.
\end{equation*} 
At the first level of construction we choose $M$ annuli inside the unit disk with the outer radius $R$ and inner radius $\bar{r}_1$ and denote their family with $\mathbb{A}_1$. As the first mapping of our construction we set 
\begin{equation*}
f_{1}(z)= \left\{
  \begin{array}{l l}
    z & \quad \text{ if $z\notin \bigcup_{B_o\in \mathbb{A}_1} B_o$}\\
   \psi_{A}(z) & \quad \text{ if $z\in B_o$ with $B_o \in \mathbb{A}_1$}.
  \end{array} \right.
\end{equation*} 
Then at the second level of construction we choose again $M$ annuli with the outer radius $\bar{R}_2$ and the inner radius $\bar{r}_2$ inside each of the inner balls of the previous generation annuli. Thus we will have $M^2$ annuli in the family $\mathbb{A}_2$ at this level and define
\begin{equation*}
\bar{f}_{2}(z)= \left\{
  \begin{array}{l l}
    z & \quad \text{ if $z\notin \bigcup_{B_o\in \mathbb{A}_2} B_o$}\\
   \psi_{A}(z) & \quad \text{ if $z\in B_o$ with $B_o \in \mathbb{A}_2$}.
  \end{array} \right.
\end{equation*} 
and 
\begin{equation*}
f_2=f_1 \circ \bar{f}_2(z).
\end{equation*}
Since $\bar{f}_2$ and $f_1$ are conformal  outside of the annuli of the corresponding level we see that the distortion of the mapping $f_2$ is non-trivial only inside the annuli from families $\mathbb{A}_1$ and $\mathbb{A}_2$. Furthermore, inside these annuli the distortion is completely determined by the  corresponding level building block \eqref{rakennuspala2}. \\
\\
We continue the construction as above. So at each level $n$ we choose $M$ new annuli from inside the previous level inner balls resulting in $M^n$ new annuli with the inner radius $\bar{r}_n$ and the outer radius $\bar{R}_n$. Denote again the family of these annuli by $\mathbb{A}_n$ and define
\begin{equation*}
\bar{f}_{n}(z)= \left\{
  \begin{array}{l l}
    z & \quad \text{ if $z\notin \bigcup_{B_o\in \mathbb{A}_n} B_o$}\\
   \psi_{A}(z) & \quad \text{ if $z\in B_o$ with $B_o \in \mathbb{A}_n$}
  \end{array} \right.
\end{equation*} 
and 
\begin{equation*}
f_n=f_{n-1} \circ \bar{f}_n(z).
\end{equation*}
Note that with a similar reasoning as in the second step the distortion is non-trivial only inside the annuli $\mathbb{A}_m$, where $m \in \{1,2,...,n\}$ and completely determined by the corresponding building block \eqref{rakennuspala2}. \\
\\
Since  radii $\bar{R}_n \to 0$ the family $f_n$ forms a Cauchy-sequence and converges uniformly towards some homeomorphism $f$. By the compactness of homeomorphisms with exponentially integrable distortion, see, for example, \cite{IM} Theorem 8.14.1 , $f$ has $p$-exponentially integrable distortion if we can show that each $f_n$ satisfies $$\int_D e^{pK_{f_n}(z)} dz < H$$
for some fixed $H$. To show this we note that the distortion of any $f_n$ is of form
\begin{equation*}
K_{f_n}(z)= \left\{
  \begin{array}{l l}
    K_{\psi_{A_{m,i}}}(z) & \quad \text{ if $z\in A_{m,i}$ for some $i\in \{1,2,...,M^m\}$ and $m\in \{1,2,...,n\}$}\\
    1 & \quad  \text{ otherwise} 
  \end{array} \right.
\end{equation*}
and hence we need to estimate the sum 
\begin{equation*}
\sum_{m=1}^{n} \int_{\cup A_{m,i}} e^{pK_{\psi_{A_{m,i}}}(z)} dz.
\end{equation*}
Since all the annuli $A_{m,i}$ are identical from the point of view of the integral we can write the sum as 
\begin{equation*}
 \sum_{m=1}^{n}  M^m \int_{B(0,\bar{R}_m)\setminus B(0,\bar{r}_m)} e^{(2-\epsilon) \log \left( \frac{1}{|z|} \right)} dz
\end{equation*}
where we have used \eqref{KKKK}. Hence straight calculation gives an estimate 
\begin{equation*}
\int_D e^{pK_{f_n}(z)} dz  \leq \pi + c\sum_{m=1}^{n}  M^m \bar{R}_{m}^{\epsilon}.
\end{equation*}
We remind that $M^m= \frac{1}{r^{sm}}$ while $$\bar{R}_{m}^{\epsilon} \leq e^{- \epsilon \left( \frac{1}{r} \right)^{\frac{(2-s)(m-1)}{p}}},$$
so thinking of $\frac{1}{r^m}$ as a variable we see that $M^m$ grows as polynomial while $\bar{R}_m$ decays like exponential. Thus we see that the the sum converges when $n \to \infty$, which ensures that the  limit map $f$ has $p$-exponentially integrable distortion. \\
\\
Let us then denote by $\mathbb{B}_n$ the set of all inner balls of the leven $n$ annuli $\mathbb{A}_n$ and set $$E=\cap_{n=1}^{\infty} \mathbb{B}_n.$$
Since each level $n$ inner ball has radius $\bar{r}_n$ and there's $M^n$ of them it is easy to see that the set $E$ has measure zero with respect to the gauge function 
$$h(t)=\frac{1}{\log^{\frac{ps}{2-s-\bar{\epsilon}}}\left( \frac{1}{t} \right)},$$
where $\bar{\epsilon}>0$. \\
\\
On the other hand each level $n$ inner ball $B_{n,i}$ is mapped under $f$ as union of $n$ linear maps imposed by \eqref{rakennuspala2}. To estimate the radius of the image ball $f(B_{n,i})$ we show that the linear stretching part of $\psi_{A_{n,i}}$  satisfies
\begin{equation}\label{LASKU5} 
\frac{\bar{R}_n \log^{\frac{p}{2-\epsilon}} \left( \frac{1}{\bar{R}_n} \right)}{\bar{r}_n \log^{\frac{p}{2-\epsilon}} \left( \frac{1}{\bar{r}_n} \right)} \bar{r}_n \geq   r^{1+\tilde{\epsilon}} \bar{r}_{n-1},
\end{equation}
where $\tilde{\epsilon}$ can be made arbitrary small by choosing $\epsilon $ small, at each level $n$. This  shows by induction that 
$$\diam f(B_{n,i}) \geq 2r^{(1+\tilde{\epsilon})n}.$$
And since we can position the level $n$ inner balls $B_{n,i}$ essentially freely we can  find a self similar set, which is built using $M$ similarity mappings with the similarity constant $ r^{1+\tilde{\epsilon}}$ and satisfies the open set condition, such that  at each level the images $f(B_{n,i})$ of the  inner balls cover the  balls of this self similar set. Thus the dimension of the set $f(E)$ is bigger or equal to the dimension of the self similar set, which is  $$\frac{s}{1+\tilde{\epsilon}}.$$ Thus, if we show \eqref{LASKU5} we  have found mapping $f$ which has $p$-exponentially integrable distortion and a set $E$ that prove the sharpness of  Theorem \ref{Korollaari}. \\
\\
To verify \eqref{LASKU5} we note that
\begin{equation*}
\begin{split}
\frac{\bar{R}_n \log^{\frac{p}{2-\epsilon}} \left( \frac{1}{\bar{R}_n} \right)}{\bar{r}_n \log^{\frac{p}{2-\epsilon}} \left( \frac{1}{\bar{r}_n} \right)} \bar{r}_n
&= \frac{R \log^{\frac{p}{2-\epsilon}} \left( \frac{1}{\bar{R}_n} \right)}{\log^{\frac{p}{2-\epsilon}} \left( \frac{1}{\bar{r}_n} \right)} \bar{r}_{n-1}\\
& \geq   \frac{R \log^{\frac{p}{2-\epsilon}} \left( \frac{1}{\bar{r}_{n-1}} \right)}{\log^{\frac{p}{2-\epsilon}} \left( \frac{1}{\bar{r}_n} \right)} \bar{r}_{n-1} \\
&=  r^{\frac{s}{2}+\frac{2-s}{2-\epsilon}} \bar{r}_{n-1} =   r^{1+\tilde{\epsilon}} \bar{r}_{n-1},
\end{split}
\end{equation*}
which finishes the proof.

\end{document}